\documentclass{ifacconf}
\usepackage{graphicx}      
\usepackage{natbib}        

\usepackage{afterpage}
\usepackage{epsfig} 
\usepackage{times} 
\usepackage{amsmath} 
\usepackage{amssymb}  
\usepackage{amsfonts}
\usepackage{mathptmx} 
\usepackage{tabularx}
\usepackage{epstopdf}
\usepackage{float}
\usepackage{epic,color}
\usepackage{mathrsfs}
\usepackage{dsfont,MnSymbol}
\usepackage{algorithm}
\usepackage{multirow} 
\usepackage[T1]{fontenc}

\newcounter{rmnum}
\newenvironment{romannum}{\begin{list}{{\upshape (\roman{rmnum})}}{\usecounter{rmnum}
			\setlength{\leftmargin}{22pt}
			\setlength{\rightmargin}{8pt}
			\setlength{\itemsep}{2pt}
			\setlength{\itemindent}{-1pt}
	}}{\end{list}}

\newcounter{anum}

\def\IEEEQEDclosed{\mbox{\rule[0pt]{1.3ex}{1.3ex}}}

\def\qed{\ifmmode\IEEEQEDclosed\else{\unskip\nobreak\hfil
		\penalty50\hskip1em\null\nobreak\hfil\IEEEQEDclosed
		\parfillskip=0pt\finalhyphendemerits=0\endgraf}\fi}

\def\qed{\hspace*{\fill}~\IEEEQED\par\endtrivlist\unskip}

\def\Re{\mathbb{R}}

\def\Sec#1{Sec.~\ref{#1}}

\def\Appendix#1{Appendix~\ref{#1}}
\def\notes#1{\marginpar{\tiny #1}\typeout{Notes!
Notes!
Notes!
}}
\renewcommand{\notes}[1]{\typeout{notes!}}
\def\FRAC#1#2#3{\genfrac{}{}{}{#1}{#2}{#3}}
\def\half{{\mathchoice{\FRAC{1}{1}{2}}%
{\FRAC{2}{1}{2}}%
{\FRAC{3}{1}{2}}%
{\FRAC{4}{1}{2}}}}

\def\Re{\field{R}}

\def\Sec#1{Sec.~\ref{#1}}

\def\clB{{\cal B}}

\def\clL{{\cal L}}
\def\clP{{\cal P}}
\def\clZ{{\cal Z}}

\def\Sec#1{Sec~\ref{#1}}

\def\E{{\sf E}}

\def\clC{{\cal C}}

\def\Sec#1{Sec.~\ref{#1}}

\def\IEEEQEDclosed{\mbox{\rule[0pt]{1.3ex}{1.3ex}}}
\def\qed{\nobreak\hfill\IEEEQEDclosed}

\def\S{O}

\def\clZ{{\cal Z}}

\newtheorem{theorem}{Theorem}

\newtheorem{definition}{Definition}

\newtheorem{remark}{Remark}
\newtheorem{proposition}{Proposition}

\def\beq{\begin{eqnarray}} 
\def\bc{\begin{center}} 
\def\be{\begin{enumerate}}
\def\bi{\begin{itemize}} 
\def\bs{\begin{small}}
\def\bS{\begin{slide}}
\def\ec{\end{center}} 
\def\ee{\end{enumerate}}
\def\ei{\end{itemize}}
\def\es{\end{small}}
\def\eS{\end{slide}}
\def\eeq{\end{eqnarray}}

\newcommand{\newP}[1]{\medskip\noindent{\bf #1:}}

\newcommand{\ud}{\,\mathrm{d}}

\def\Re{\mathbb{R}}
\def\E{{\sf E}}

\def\clY{{\cal Y}}

\def\Sec#1{Sec.~\ref{#1}}

\def\clB{{\cal B}}

\def\clL{{\cal L}}
\def\clP{{\cal P}}
\def\clZ{{\cal Z}}

\renewcommand{\Re}{\mathbb{R}}

\def\FRAC#1#2#3{\genfrac{}{}{}{#1}{#2}{#3}}

\def\sP{{\sf P}}

\def\Nsp{{\sf N}} 
\def\Rsp{{\sf R}} 
\def\ones{{\sf 1}}

\def\clA{{\cal A}}
\def\clC{{\cal C}}

\def\clL{{\cal L}}
\def\clO{{\cal O}}
\def\clP{{\cal P}}
\def\clU{{\cal U}}
\def\clY{{\cal Y}}
\def\clZ{{\cal Z}}
\def\E{{\sf E}}
\def\hE{\tilde{\sf E}}
\def\bP{{\sf P}}

\def\bS{\mathbb{S}}
\def\dv{\operatorname{diag}}
\def\sp{\operatorname{span}}
\def\bpi{\tilde{\pi}}

\newlength{\noteWidth}
\long\def\notes#1{\ifinner
	{\tiny #1}
	\else
	\marginpar{\parbox[t]{\noteWidth}{\raggedright\tiny #1}}
	\fi}

\long\def\notes#1{}
\begin{document}
	\begin{frontmatter}
		
		\title{A Dual Characterization of Observability for Stochastic Systems\thanksref{footnoteinfo}} 
		
		\thanks[footnoteinfo]{Financial support from the 
			NSF grant 1761622 and the ARO grant W911NF1810334 is gratefully acknowledged. }
		
		\author[First]{Jin Won Kim} 
		\author[First]{Prashant G. Mehta} 
		
\address[First]{Coordinated Science Laboratory and Department of Mechanical Science and Engineering, University of Illinois at Urbana-Champaign, Urbana, IL 61801 USA (e-mail: jkim684@illinois.edu; mehtapg@illinois.edu)}

		\begin{abstract}                
			This paper is concerned with the definition and characterization of the observability for a
			continuous-time hidden Markov model where the state evolves as a
			continuous-time Markov process on a compact state space and the observation process is modeled as nonlinear function of the state corrupted by a Gaussian measurement noise.
			The main technical tool is based on the recently
			discovered duality relationship between minimum variance estimation
			and stochastic optimal control: The observability is defined as a dual
			of the controllability for a certain backward stochastic differential
			equation.  Based on the dual formulation, a
                        test for observability is presented and
                        related to 
                        literature. 
			The proposed duality-based framework allows
                        one to easily relate and compare the linear
                        and the nonlinear systems.  A side-by-side
                        summary of this relationship is given in a
                        tabular form (Table~1).
			

		\end{abstract}
		
		
	\end{frontmatter}

\section{Introduction}

This paper is concerned with the definition of
observability for a partially observed pair of
continuous-time stochastic processes $(X,Z)$ where the state $X$ is a
Markov process and the observation $Z$ is a nonlinear function of the
state corrupted by the Gaussian measurement noise.  The precise
mathematical model appears in the main
body of the paper.


In deterministic linear time-invariant (LTI) settings, observability (more generally
detectability) and
its dual relationship to the controllability are
foundational concepts in linear systems theory;~\citep{kalman1960general}.  It is an important property that a
model must satisfy to construct asymptotically stable
observers~\citep{kailath2000linear}.  
For a partially observed stochastic LTI 
system, the detectability property of its deterministic counterpart
is necessary and also sufficient (under mild additional conditions) to deduce 
results on
asymptotic stability of the optimal (Kalman) filter;~\citep{ocone1996}.

Generalization of these concepts to nonlinear deterministic and
stochastic systems has been an area of historical and current research
interest~\citep{hermann1977nonlinear,krener1985nonlinear,moraal1995observer,liu2011stochastic}.
In settings more general than this paper, the fundamental definition
of observability is due to \cite{van2009observability,van2009uniform}.
The definition is used to establish results on asymptotic stability of
the nonlinear filter~\citep{chigansky2009intrinsic,van2010nonlinear}.
Certain extensions of van Handel's observability definition appear in 
recent papers by~\cite{mcdonald2018stability,mcdonald2019cdc}.

In this paper, we utilize the recently discovered duality relationship
between minimum variance estimation and stochastic optimal
control (see~\cite{kim2019bsdeduality}) to define observability as a
dual to the controllability.  The latter property is somewhat
`natural' because it bears close resemblance to the definition of
controllability in deterministic LTI settings.  The definition of
observability is obtained by using duality. 
In finite state-space settings, certain Kalman-type rank conditions are derived
to verify the observability property.  These conditions are shown to be identical to the
ones reported by~\cite{van2009observability} but derived
here using alternate means.  

Given the close similarity of the finite state-space hidden Markov
model and the deterministic LTI model, these conditions allow one to
compare and contrast the differences between the two.
This is important given many attempts and successes over the years to apply methods and tools from linear systems theory to study finite-dimensional Markov and hidden Markov models;
cf.,~\citep{kotsalis2008balanced,Deng_cdc10,Brockett_cdc09,chebusmey18a}.

\begin{table*}[t]
	\centering
	\renewcommand{\arraystretch}{1.7}
	\small
	\caption{Comparison of the controllability--observability duality for linear and nonlinear systems}
	\label{tb:comparison}
	\begin{tabular}{m{0.13\textwidth}m{0.4\textwidth}m{0.4\textwidth}}
		& {\bf Linear-deterministic case} (\Sec{sec:background}) & {\bf Nonlinear-stochastic case} (\Sec{sec:main}) \\ \hline
		Signal space & $\clU =  L^2([0,T];\Re^m)$\vspace{3pt} \newline  $\langle u,v\rangle_\clU = \displaystyle\int_0^T u_t^\top v_t\ud t$ \vspace{2pt} & $\clU =  L^2_\clZ([0,T];\Re^m)
		$\vspace{3pt} \newline $\langle U,V\rangle_\clU = \displaystyle\hE\Big(\int_0^T U_t^\top V_t\ud t\Big)$\vspace{2pt} \\ \hline
		Function space & $\clY =\Re^d$ \vspace{2pt} \newline $\langle x,y \rangle_\clY = x^\top y$ & $\clY = C(\bS)$, $\clY^\dagger = {\cal M}(\bS)$ \vspace{2pt} \newline $\langle \mu,y\rangle_\clY = \mu(y)$ \\ \hline
		Linear operator & $\clL : \clU \to \clY$ \vspace{2pt} \newline
		\phantom{$\clL$} \hspace{0.2em} $u\mapsto y_0$ by ODE~\eqref{eq:LTI-ctrl} & $\clL:\clU\times\Re\to\clY$ \vspace{2pt} \newline
		\phantom{$\clL$} \hspace{0.2em}$(U,c)\mapsto Y_0$ by BSDE~\eqref{eq:NL-ctrl} \\ \hline
		Adjoint operator & $\clL^\dagger : \clY \to \clU$ \newline
		\phantom{$\clL$} \hspace{0.15em} $x_0\mapsto z_t$ by ODE~\eqref{eq:LTI-obs} & $\clL^\dagger:\clY^\dagger\to \clU\times\Re$, \vspace{2pt} \newline \phantom{$\clL^\dagger $} \hspace{0.15em} $\bpi_0\mapsto (\bpi_t(h), \bpi_0(\ones))$ by Zakai equation~\eqref{eq:Zakai} \\ \hline
		Observability & $\Rsp(\clL) = \clY \quad  \iff \quad \Nsp(\clL^\dagger) = \{0\}$  \newline $\phantom{\Rsp(\clL) = \clY \quad}\iff$ \vspace{2pt} \newline $ He^{At}x_0^{(1)} \equiv He^{At}x_0^{(2)} \Rightarrow x_0^{(1)} = x_0^{(2)}$ \quad\eqref{eq:LTI-O3} & $\overline{\Rsp(\clL)} = \clY \quad  \iff \quad \Nsp(\clL^\dagger) = \{0\}$  \newline $\phantom{\Rsp(\clL) = \clY \quad}\iff$ \vspace{2pt} \newline $\E^\mu (h(X_t)|\clZ_t) = \E^\nu(h(X_t)|\clZ_t)
		\Rightarrow \mu = \nu$ \quad\eqref{eq:emu-enu} \\ \hline
	\end{tabular}
\end{table*}

The original contributions of our paper are as follows:
\begin{romannum}
\item The use of dual control problem to define observability and obtain 
  conditions for the same is original.  The type of duality used in
  our work was introduced for the first time in our recent paper~\citep{kim2019bsdeduality}.  As fully
  explained in the introduction to this earlier paper, our notion of
  \emph{duality} is {entirely different} 
   from prior works on duality such as:~\cite{mortensen1968,simon1970duality,FlemingMitter82,fleming1997deterministic,mitter2003,goodwin2005,todorov2008general}. 
\item Although we relate our definition of
  observability to~\cite{van2009observability}, the background and
  motivation for our work is different.  We
  essentially define observability in terms of the controllability property of the
  dual.  For nonlinear stochastic systems, such an approach is
  original and different from Van Handel's approach which is entirely probabilistic in nature.  In particular, there is no hint
  in any of the earlier works to use duality to define observability
  for nonlinear stochastic systems. (In contrast, the duality
  between controllability and observability is standard for linear systems).  
\item The upshot of our work is that we can establish parallels
  between linear and nonlinear cases (see Table~\ref{tb:comparison}).
  This is expected to be useful in several ways, e.g., to obtain
  approximation algorithms and for stability analysis of the filter.  

\end{romannum}

The remainder of this paper is organized as follows: The background on
the classical deterministic LTI model appears in \Sec{sec:background}. The nonlinear model is introduced in \Sec{sec:model} and its stochastic observability defined and discussed in \Sec{sec:main}.  The finite state case is illustrated in \Sec{sec:finite}. The conclusions appear in
\Sec{sec:conclusion}.  All the proofs are contained in the Appendix.    




\section{Background: Duality in linear systems}\label{sec:background}

In linear algebra, it is an elementary fact that the range space of a
matrix is orthogonal to the null space of its transpose.  
In functional analysis, the closed range theorem provides the
necessary generalization of this elementary fact in
infinite-dimensional settings.  The theorem~\citep[Theorem
6.5.10]{hutson2005applications} states that 
$$
\overline{\Rsp(\clL)} = \Nsp(\clL^\dagger)^\bot
$$
where $\overline{\Rsp(\clL)}$ is closure of the range space of a bounded linear operator $\clL$ and 
$\Nsp(\clL^\dagger)$ is the null space of its adjoint operator
$\clL^\dagger$.  This dual relationship is of fundamental importance
to understand the duality between controllability and observability.  
The overall procedure is as follows: 
\begin{enumerate}
\item Define the appropriate function spaces and the associated linear
  operators; and
\item Express
controllability and observability properties in terms of range space
and null space of these operators. 
\end{enumerate} 
We briefly review this well known procedure first
in the classical settings.


\newP{Function spaces} Denote ${\cal U} := L^2([0,T];\Re^m)$ to be the
Hilbert space of $\Re^m$-valued (input or output) square-integrable
signals on the time interval $[0,T]$.  The space is equipped with the
inner product $\langle u,v\rangle_{\cal U} = \int_0^T u_t^\top v_t \ud t$ for 
$u,v\in\clU$.  Denote ${\cal Y}:=\Re^d$ to be the Euclidean space
equipped with the standard inner product $\langle y_0,
x_0\rangle_{\cal Y}:=y_0^\top x_0$ for $y_0, x_0 \in\clY$.


\newP{Operators} For given matrices $A\in\Re^{d\times d}$ and
$H\in\Re^{m\times d}$ define a linear operator $\clL:\clU\to\clY$ as follows:
\[
\clL u = \int_0^T e^{A^\top t}H^\top u_t \ud t =: y_0
\]
The definition of the adjoint operator $\clL^\dagger: \clY \to \clU$
follows from the following calculation:
\begin{equation}
\langle \clL u,
x_0\rangle_{\cal Y} = y_0^\top x_0 = \int_0^T
u_t^\top H e^{At}x_0\ud t = \langle u,\clL^\dagger x_0\rangle_{\cal U}
\label{eq:dual_LTI}
\end{equation}
Therefore,
\[
\big(\clL^\dagger x_0\big)(t) = H e^{At} x_0=:z_t\quad \text{for}\;\;t\in[0,T]
\]

\newP{Controllability and observability} The operator $\clL$ defines
the map from a given input signal $u=\{u_t:0\le t\le T\}$ to the
initial condition $y_0$ for the linear system\footnote{The system~\eqref{eq:LTI-ctrl} is an example of a backward ordinary differential equation (ODE)
because the terminal condition at time $t=T$ is set (to zero in this
case).}:
\begin{equation}\label{eq:LTI-ctrl}
-\dot{y}_t = A^\top y_t + H^\top u_t,\quad y_T = 0
\end{equation}
The range space $\Rsp(\clL)$ is referred to as the {\em
  controllable subspace}.  The system~\eqref{eq:LTI-ctrl} is said to
be {\em controllable} if $\Rsp(\clL)={\cal Y}$.

The adjoint operator $\clL^\dagger$ defines the map from a given
initial condition $x_0$ to the observation signal $z=\{z_t:0\le t\le
T\}$ for the linear system
\begin{subequations}\label{eq:LTI-obs}
	\begin{align}
	\dot{x}_t &= Ax_t,\quad \text{with init. cond.}\;\; x_0 \label{eq:LTI-obs-a}\\
	z_t &= H x_t \label{eq:LTI-obs-b}
	\end{align}
\end{subequations}
The system~\eqref{eq:LTI-obs} (henceforth referred to as the linear model
$(A,H)$) is said to be {\em observable} if
$\Nsp(\clL^\dagger)=\{0\}$.  

By the closed-range theorem (or more directly by simply using~\eqref{eq:dual_LTI}), $
\Rsp(\clL) = \Nsp(\clL^\dagger)^\perp$.  Therefore, the
system~\eqref{eq:LTI-obs} is observable if and only if the 
system~\eqref{eq:LTI-ctrl} is controllable.  This is useful in the
following ways:
\begin{romannum}
\item Definition of observability: as the property of the dual
  system being controllable.  
\item Geometric interpretation of non-observability: If the controllable subspace $\Rsp(\clL) \subsetneq \Re^d$ then there
exists a non-zero vector $\tilde{x}_0\in\Nsp(\clL^\dagger)$ such that
$
y_0^\top \tilde{x}_0 = 0$ for all $y_0\in\Rsp(\clL)
$. 
The vector $\tilde{x}_0$ has an interpretation of being the
``un-observable'' direction in the following sense: For any given
$x_0\in\Re^d$,
$
He^{At}x_0 \equiv He^{At}(x_0+\tilde{x}_0)$ for all $t\in[0,T]
$.  This in turn provides an equivalent definition of observability:
The model $(A,H)$ is observable if
\begin{equation}\label{eq:LTI-O3}
He^{At}x_0^{(1)} \equiv He^{At}x_0^{(2)} \; \forall\; t\in[0,T] \text{ implies } x_0^{(1)} = x_0^{(2)}
\end{equation} 





\item Tests for observability:  
By the use of the Cayley-Hamilton theorem, 
\begin{equation}\label{eq:obs_gram}
	\Rsp(\clL) = \sp\big\{H^\top, A^\top H^\top , \ldots, (A^\top)^{d-1} H^\top \big\}
	\end{equation}
This provides a straightforward test to verify observability: 
The model $(A,H)$ is observable if the span on the right-hand side 
of~\eqref{eq:obs_gram} is $\Re^d$.  
\end{romannum}

The aim of this paper is to repeat the above program---viz., (i) the
definition of the function spaces ${\cal U}$ and ${\cal Y}$; (ii) the
definition of the linear operator $\clL$ and its adjoint
$\clL^\dagger$; (iii)
the mathematical characterization of the controllable subspace
$\Rsp(\clL)$; and (iv) its use in definition and geometric interpretation of
the observability---for a partially observed nonlinear stochastic
system.  A summary of the paper appears in the form of a
comparison between the linear-deterministic and the
nonlinear-stochastic system in Table~\ref{tb:comparison}.


\section{Problem Formulation}\label{sec:model}


\def\lsq#1#2{L_{#1}^2([0,T]\,;{#2})}




\subsection{Model \& notation} 

The nonlinear model is defined for a pair of continuous-time
stochastic processes denoted as $(X,Z)$.  The 
details of the model are as follows:   
\begin{romannum}
\item The state $X=\{X_t:\, 0\le t\le T \}$ is a Markov process that
  evolves in a compact state-space $\bS$\footnote{The results of this
    paper are expected to carry over to locally compact spaces. }. The generator of the Markov
  process $X$ is denoted as $\clA$.  
\item The observation process
$Z=\{ Z_t :\, 0\le t\le T \}$ is defined according to the following model:
\begin{equation}\label{eq:obs_model}
Z_t = \int_0^t h(X_s) \ud s + W_t
\end{equation}
where $h:\bS\rightarrow \Re^m$ is a given observation function and
$\{W_t:\,t \ge 0\}$ is an $m$-dimensional Wiener process (w.p.). It is
assumed that $W$ is independent of
$X$.
\item We refer to the above model as the nonlinear model
  $(\clA,h)$.  
\end{romannum}

\newP{Notation} We denote $\clZ_t := \sigma(\{Z_s:s\le t\})$ to be the
  $\sigma$-algebra generated by the observations up to time $t$ and $\clZ
:= \{\clZ_t:0\le t\le T\}$ is the entire filtration.

The law of $(X,Z)$ is denoted as $\bP$ with the associated expectation operator $\E$.  To emphasize the model for initial
  condition $X_0$, we use $\bP^\mu$ to denote the law of $(X,Z)$ when
  the initial condition $X_0\sim \mu$.  

For the state-space $\bS$, we let $\clB(\bS)$ denote the
  Borel $\sigma$-algebra on $\bS$; ${\cal M}(\bS)$ is the
  vector space
  of (signed) Radon (bounded and regular) measures on $\clB(\bS)$; and
  $\clP(\bS)\subset {\cal M}(\bS)$ is the set of probability
  measures.  $C(\bS)$ is used to denote the dual space of ${\cal
    M}(\bS)$, which is identified by continuous functions on
  $\bS$~\citep[Ch. 21]{treves1967topological}. Throughout this paper, we will use the notation:
  $$\mu(f):=\int_{\bS} f(x) \mu(\ud x)$$
  to denote the integral of a measurable function $f$ with respect to the measure $\mu$. It is the natural duality paring between ${\cal M}(\bS)$ and $C(\bS)$.

%
%
%

\subsection{Preliminaries} 

The main concern of this paper is to define and characterize
observability for the nonlinear model $(\clA,h)$. In his paper,
\cite{van2009observability} proposes the following probabilistic
definition\footnote{The definition in
\cite{van2009observability} applies to a more general
class of stochastic processes $(X,Z)$ whereby the independent increment of the
measurement noise may not be of the additive Gaussian form (as assumed
here).} of observability for stochastic processes $(X,Z)$:

\medskip

\begin{definition}
Suppose $X$ is a Markov process defined on a compact set $\bS$ and $Z$ is defined according to
model~\eqref{eq:obs_model}.
Suppose ${\sf P}^\mu$ and ${\sf P}^\nu$
are two laws of the process $(X,Z)$ with initial measure $X_0\sim\mu$
and $X_0\sim\nu$, respectively.  The model is said to be {\em P-observable} if   
\begin{equation}\label{eq:VH-obs-defn-eqn}
	\bP^{\mu}\big|_{\clZ_T} = \bP^{\nu}\big|_{\clZ_T} \quad
        \Rightarrow \quad \mu = \nu
\end{equation}
where $\bP^{\mu}\big|_{\clZ_T}$ denotes the restriction of the
probability measure $\bP^{\mu}$ to the $\sigma$-algebra $\clZ_T$.  
\end{definition} 

\medskip
  
%
%

Before presenting the main result, it is
useful to review some concepts from the theory of nonlinear
filtering~\citep[Ch. 5]{xiong2008introduction}:

\newP{Change of measure} 
Given $\bP$, define a new measure $\tilde{\bP}$ according to the
Radon-Nikodyn derivative
\begin{equation*}
\frac{\ud \tilde{\bP}}{\ud \bP}(\omega) := \exp\Big(-\int_0^T
h^\top(X_t)\ud Z_t + \half\int_0^T |h(X_t)|^2 \ud t\Big)
\end{equation*}
By the Girsanov theorem, $Z$ is a $\tilde{\bP}$ Wiener process. 
For a given function $f$, the un-normalized filter is defined by
\begin{equation*}\label{eq:unnormalized-filter}
\sigma_t(f) := \hE(D_tf(X_t)|\clZ_t)
\end{equation*}
where $\tilde{\E}(\cdot)$ denotes the expectation operator with respect to the
new measure $\tilde{\bP}$ and 
$$
D_t
= \exp\Big(\int_0^t h^\top(X_s)\ud Z_s - \half\int_0^t |h(X_s)|^2 \ud s\Big)
$$ 
The un-normalized filter $\sigma_t(f)$ solves the Zakai equation of
nonlinear filtering.  The nonlinear filter is given by
\begin{equation}\label{eq:normalize-filter}
\pi_t(f):=\E(f(X_t)|\clZ_t) = \frac{\sigma_t(f)}{\sigma_t(\ones)}
\end{equation}
where $\ones(x)\equiv 1 \; \forall \; x\in\bS$ denotes the unit
constant function. As before, we use superscript
(e.g.,~$\sigma_t^\mu$, $\pi_t^\mu$) to emphasize dependence on the initial
measure ($\mu$) of $X_0$.

\medskip

\section{Main result: Stochastic observability}\label{sec:main}


\subsection{Function spaces}
\label{sec:FS_nl}

In nonlinear settings, the signal space $\clU = L^2_\clZ([0,T];\Re^m)$
is the Hilbert space of $\Re^m$-valued stochastic processes
on $[0,T]$.  The subscript $\clZ$ denotes the fact that the signals
are (forward) adapted to the filtration $\clZ$.  
The space is equipped with the inner product
$$
\langle U, V\rangle_\clU := \hE\;\Big(\int_0^T U_t^\top V_t \ud t\Big)
$$ 
The expectation $\hE$ is with respect to the measure
$\tilde{\bP}$. For the proof that $\clU$ is a Hilbert space with
respect to this inner product~\citep[p. 99]{le2016brownian}.  



The space $\clY = C(\bS)$ and its dual $\clY^\dagger={\cal M}(\bS)$.
For a function $y\in C(\bS)$ and a measure $\mu\in{\cal
  M}(\bS)$, the 
dual pairing is as follows:
$$
\langle \mu,f\rangle_\clY = \mu(f) = \int_{\bS} f(x) \mu(\ud x) 
$$

A side-by-side comparison of the signal space and the function space
for the linear and nonlinear cases appears as first two rows in
Table~1.  

\subsection{Controllability}

Parallel to the linear case, we define controllable subspace as the range space of
a bounded linear operator.  For this purpose, we introduce the
following 
backward stochastic differential equation (BSDE):
\begin{align}
-\ud Y_t(x) &= \big(\clA Y_t(x) + h^\top(x)(U_t + V_t(x))\big)\ud t - V_t^\top(x)\ud Z_t \nonumber\\
Y_T(x) & = c \ones(x)  \;\;\forall \; x\in\bS 
      \label{eq:NL-ctrl}
\end{align}
where $c\in\Re$ and the input signal $U\in\clU$.  The
solution $(Y,V):=\{(Y_t(x),V_t(x))\,:\, t \in [0,T],\;x\in \bS\}$ of
the BSDE is (forward) adapted to the filtration $\clZ$.  For the purposes of this paper,
well-posedness (existence, uniqueness and regularity) of the solution $(Y,V)\in L^2_{\clZ}([0,T];C(\bS))\times
L^2_{\clZ}([0,T];C(\bS))$ is assumed;
cf.,~\citep{ma1997adapted}.  The BSDE is the nonlinear counterpart of
the backward ode~\eqref{eq:LTI-ctrl} in the LTI setting.
The justification for considering the BSDE~\eqref{eq:NL-ctrl} appears in~\Appendix{apdx:bsdeduality} where our prior work~\citep{kim2019bsdeduality} on the topic of duality is briefly reviewed.  

The bounded linear operator\footnote{The bounded-ness property is based on the
  well-posedness of the solution of the BSDE~\eqref{eq:NL-ctrl}.}
$\clL:\clU\times\Re\to\clY$ is defined through the solution of the BSDE~\eqref{eq:NL-ctrl} as follows:
\begin{equation}\label{eq:L_nl_defn}
\clL(U,\;c) = Y_0
\end{equation}
and its range space $
\Rsp(\clL) = \{Y_0\in \clY:U\in \clU,c\in\Re\}
$ is referred to as the {\em controllable space}.   The
BSDE~\eqref{eq:NL-ctrl} is said to be {\em controllable} if
$\Rsp(\clL)$ is dense in $\clY$.  

In finite state-space settings, when the state space $\mathbb{S}$ is
of cardinality $d$, $\Rsp(\clL)$ is a subspace of $\Re^d$.
Therefore, in this setting, the system is {\em controllable} if
$\Rsp(\clL)=\Re^d$.  

Duality is used to propose an indirect definition
of observability as follows:



\medskip

\begin{definition}\label{def:indirect}
The nonlinear model $(\clA,h)$ is said to be
{\em observable} if 
\begin{equation}\label{eq:O1}
\Rsp(\clL) \text{ is dense in }\clY \tag{O1} 
\end{equation}  
\end{definition}



\subsection{Observability}

We develop a more direct definition of observability by considering the dual operator. In the Prop.~\ref{thm:duality} (stated below), it is shown that the adjoint to the BSDE~\eqref{eq:NL-ctrl} is the \emph{Zakai equation}:
\begin{equation}\label{eq:Zakai}
\bpi_t(f) = \bpi_0(f) + \int_0^t \bpi_s(\clA f)\ud s + \int_0^t \bpi_s(h^\top f)\ud Z_s\quad \forall f \in \clY
\end{equation}
where the initial condition $\bpi_0\in{\cal M}(\bS)$ is
given\footnote{In nonlinear filtering, the Zakai equation is
  considered with initial measure $\bpi_0\in{\cal P}(\bS)$.  In this
  paper, the
  initial measure is allowed to be a signed measure.}. For a
given function $f\in\clY$, the
solution of the Zakai equation~\eqref{eq:Zakai} is denoted as
$\bpi(f):=\{\bpi_t(f):0\le t\le T\}$.  In finite
state-space settings, the Zakai equation is simply a linear SDE on
$\Re^d$ with initial measure $\bpi_0\in\Re^d$.

The following proposition is proved in~\Appendix{apdx:proof-thm1}:

\medskip

\begin{proposition}\label{thm:duality}
Consider the linear operator~\eqref{eq:L_nl_defn}.  Its adjoint
$\clL^\dagger: \clY^\dagger \to \clU \times \Re$ is given by
\[
\clL^\dagger \bpi_0 = (\bpi(h),\;\bpi_0(\ones))
\]
where $\bpi(h)=\{\bpi_t(h):0\le t\le T\}$ is the solution of the Zakai equation~\eqref{eq:Zakai}
with $f=h$ and the initial measure $\bpi_0\in \clY^\dagger$.   
\end{proposition}

\medskip

For the purposes of defining observability, the adjoint's null space $
\Nsp(\clL^\dagger) = \{\bpi_0\in\clY^\dagger: \bpi(h) = 0, \; \bpi_0(\ones) = 0\}
$ is of interest.  In the finite
state-space settings, $\Nsp(\clL^\dagger)$ is a subspace of $\Re^d$.

The dual of definition~\eqref{eq:O1} is as follows:

\medskip

\begin{definition}\label{def:direct}
The nonlinear model $(\clA,h)$ is said
to be {\em observable} if
\begin{equation}\label{eq:O2}
\Nsp(\clL^\dagger) = \{0\} \tag{O2}
\end{equation}
\end{definition}

\medskip

The two definitions~\eqref{eq:O1} and~\eqref{eq:O2} are equivalent: By the closed range theorem $
\overline{\Rsp(\clL)} = \Nsp(\clL^\dagger)^\perp$.  
If the controllable subspace $\overline{\Rsp(\clL)} \subsetneq \clY$
then there exists a non-zero measure ${\bpi}_0\in\Nsp(\clL^\dagger)$ such that
$
{\bpi}_0(Y_0) = 0$ for all $Y_0\in\overline{\Rsp(\clL)}
$. 
The measure $\bpi_0$ has an interpretation of being the
\emph{un-observable measure} in the following sense: For given
$\mu\in{\cal P}(\bS)$ being a ``true'' distribution of $X_0$, suppose $\epsilon\neq 0$ is chosen such
that $\nu = \mu + \epsilon \bpi_0 \in {\cal P}(\bS)$. Then owing to the linearity of~\eqref{eq:Zakai},
$$\sigma_t^\mu(h) = \sigma_t^\nu(h)\quad t\text{-a.e.}\;\;\tilde{\sP}^\mu\text{-a.s.}$$
As will be justified more fully in the proof of
Theorem~\ref{prop:equivalent-definition}, this leads to the third
equivalent definition of observability:

\begin{definition}\label{def:preocess}
The nonlinear model $(\clA,h)$ is said
to be {\em observable} if
\begin{equation}\label{eq:emu-enu}
	\pi_t^\mu(h) = \pi_t^\nu(h) \quad t\text{-a.e.}\;\; \sP^\mu\text{-a.s.}\quad
        \Rightarrow \quad \mu = \nu \tag{O3}
\end{equation}
\end{definition}

It is noted that~\eqref{eq:emu-enu} is the stochastic analog
of~\eqref{eq:LTI-O3}. 

The proof of the following theorem appears in the Appendix~\ref{apdx:unobservable-direction}.

\begin{theorem}[Observability]\label{prop:equivalent-definition}
The three conditions:~\eqref{eq:O1},~\eqref{eq:O2}, and \eqref{eq:emu-enu} are equivalent.  
\end{theorem}

\subsection{Test for observability}



The following theorem provides an explicit characterization of
the controllable space.  Its proof appears in the Appendix~\ref{apdx:proof-rank}.  

\medskip

\begin{theorem}\label{thm:rank}
Consider the linear operator~\eqref{eq:L_nl_defn}.  Its range space
${\Rsp(\clL)}$ is the smallest such subspace $\clC\subset\clY$ that satisfies the
following two properties:
\begin{romannum}
	\item The constant function $\ones\in \clC$;
	\item If $g\in\clC$ then $\clA g \in \clC$ and $g \cdot h
          \in\clC$. ($g \cdot h$ is the Hadamard (element-wise)
          product of functions $g$ and $h$)\footnote{For a
            vector-valued function $h(x) = [h_1(x),\ldots,h_m(x)]$,
            $g\cdot h\in\clC$ means $g\cdot h_i \in \clC$ for each
            $i=1,\ldots,m$.  The Hadamard product is simply the product of functions, i.e., $(g\cdot h_i)(x)=g(x)h_i(x)$
            for all $x\in\bS$.}.
\end{romannum}
\end{theorem}

\medskip

\subsection{Relationship to P-observability}




For the particular (additive Gaussian noise) form of the observation
model~\eqref{eq:obs_model}, there is a formula, due to~\cite[Theorem
3.1]{clark1999relative}, for the relative entropy
between $\bP^{\mu}_{\clZ_T}$ and $\bP^{\nu}_{\clZ_T}$:
\begin{equation*}\label{eq:relative-entropy}
{\sf D} (\bP^{\mu}_{\clZ_T} \| \bP^{\nu}_{\clZ_T}) =
\E^{\mu} \int_0^T |\pi_t^{\mu}(h) - \pi_t^{\nu}(h)|^2 \ud t 
\end{equation*}
Combined with Theorem~\ref{prop:equivalent-definition}, a
straightforward corollary is the following
proposition:

\medskip

\begin{proposition}\label{thm:equivalency}
	Consider the observation model
	of the form~\eqref{eq:obs_model}.  The model $(\clA, h)$ is
	observable (according to one of the equivalent definitions~2, 3, or 4)
	if and only if it is P-observable (definition~1).
\end{proposition}


\section{Finite state-space case}\label{sec:finite}



The results of this paper are next illustrated for the finite
state-space case.  
The following notation is adopted:
\begin{romannum}
	\item Without loss of generality, it is convenient to consider the
	state space $\bS= \{e_1,\ldots,e_d\}$ defined by the
	canonical basis in $\Re^d$.
	\item The function space
	$C(\bS)$ is identified with $\Re^d$: any function $\tilde{f}:
	\bS\to \Re$ is determined by its values at the basis vectors
	$\{e_i\}$.  We denote these values as a column vector $f\in\Re^d$
	and express $\tilde{f}(x)=f^\top x$ for $x\in \bS$.  In this section, with a slight abuse of notation, we will drop the tilde
	to simply write $f(x)=f^\top x$.  
	\item The set ${\cal P}(\bS)$ is the probability simplex in $\Re^d$.  For a measure $\mu \in {\cal M}(\bS)$, the integral
	$\mu(f)=\sum_i \mu(e_i) f(e_i) = \mu^\top f$ is a dot product. 
	\item The observation function $h(x) =
	H^\top x$, where $H\in\Re^{d\times m}$. 
	\item The generator $\clA$ of the Markov process is identified with a
	row-stochastic rate matrix $A\in\Re^{d\times d}$ which acts on
	functions (elements of $\Re^d$) through right-multiplication: $A: f
	\mapsto A f$.
 	\item Given the identification of the generator $\clA$ with
          the matrix $A$ and the
	observation function $h(x)$ with the matrix $H$, we refer to
        the finite-state Markov chain
	as the nonlinear model $(A,H)$.
\end{romannum}

Using this notation, the BSDE~\eqref{eq:NL-ctrl} is expressed as follows:
\begin{equation*}\label{eq:dual-ctrl-finite}
-\ud Y_t = \Big(AY_t +HU_t+ \dv^\dagger(HV_t^\top)\Big)\ud t - V_t \ud Z_t,\quad Y_T = c\ones
\end{equation*}
where $\ones$ is a vector of ones in $\Re^d$ and
$\dv^\dagger(HV_t^\top)$ is the vector of the diagonal elements of the matrix $HV^\top$. The solution pair is $(Y,V) \in L^2_\clZ([0,T];\Re^d)\times L^2_\clZ([0,T];\Re^{d\times m})$.


The controllable space $\Rsp(\clL)$ is a subspace in $\Re^d$:
\begin{align}
\Rsp(\clL) = \sp\big\{\ones, &\,  H, \,  AH, \,  A^2H, \,  A^3H, \, \ldots, \label{eq:obs_gram_nl}\\
&H\cdot H, \,  A(H \cdot H), \,  H\cdot (AH), \,  A^2(H\cdot H),\ldots, \nonumber\\
&H\cdot (H\cdot H), \,  (AH)\cdot (H\cdot H), \,  H\cdot A(H \cdot H), \, \ldots \big\} \nonumber
\end{align}
where the dot denotes the element-wise product.
The nonlinear model $(A,H)$ is observable if the vectors in the
righthand-side of~\eqref{eq:obs_gram_nl} span $\Re^d$.  This provides
a test for verifying observability of the nonlinear model.  

We next compare the above test with the observability
test~\eqref{eq:obs_gram} for the linear
model $(A,H)$.  It is clear that if the linear model $(A,H)$ is
observable (in the sense of~\eqref{eq:obs_gram}) then the nonlinear
model $(A,H)$ is also observable.  However, the latter property is in general much
weaker than the observability in linear systems theory.  This is shown in the following
proposition whose proof appears in~\Appendix{apdx:injective-h}.









 \medskip

\begin{proposition}[A sufficient condition]\label{prop:sufficient}
Consider the nonlinear model $(A,H)$ for the finite state-space.
Then $(A,H)$ is observable if $h(x)=H^\top x$ is an injective map from
$\bS$ into $\Re^m$. (The map is injective if and only if $H_i \neq H_j$ for all $i\neq
j$ where $H_i$ is the $i^\text{th}$ row of the $d\times m$ matrix $H$).    
If $A=0$ then the injective property of the function $h$ is
also necessary for observability.    
\end{proposition}

\medskip

\begin{remark}\label{rm:rank-condition-remarks}
For the finite state nonlinear model $(A,H)$, the test for
observability first appeared in~\cite [Lemma~9]{van2009observability}.
The test was obtained by explicitly calculating the probability of the event
$[h(X_{t_1})=n_1,  h(X_{t_2})=n_2, \ldots, h(X_{t_k})=n_k]$ and
applying~\eqref{eq:VH-obs-defn-eqn}; cf.,~\citep[Lemma 8]{van2009observability}.    
For a general class of linear BSDE-s, the controllable subspace is identically defined by~\cite[Lemma 3.2]{peng1994backward}.  However, its use in the study of observability appears to be new.  
\end{remark}

 \section{Conclusion and Future Directions}\label{sec:conclusion}

In this paper, the duality introduced in our recent
work~\citep{kim2019bsdeduality} is used for the purposes of defining
and characterizing observability of nonlinear stochastic systems.  The
main idea is to define observability as a dual of the controllability
of a certain BSDE~\eqref{eq:NL-ctrl}.   Based on the dual formulation, a
                        test for observability is presented and
                        related to 
                        literature. 
			The proposed duality-based framework allows
                        one to easily relate and compare the linear
                        and the nonlinear systems.  A side-by-side
                        summary of this relationship is given in a
                        tabular form (Table~1).

The methodology of this and our earlier duality paper is currently being used to investigate nonlinear filter stability;
and to develop new control-based algorithms for approximating the
nonlinear filter.  This is the subject of continuing research.    

\bibliographystyle{ifacconf}
\bibliography{duality,backward_sde,filter-stability-observability}
\appendix

\section{Duality between estimation and control}\label{apdx:bsdeduality}

This section includes a brief review of the duality between nonlinear
filtering and stochastic optimal control introduced in our recent
paper~\citep{kim2019bsdeduality}.


\newP{Dual optimal control problem}  
\begin{align}
&\mathop{\text{Min}}_{U\in\,{\cal U}}\ {\sf J}(U) =
\E\Big(\half|Y_0(X_0)-\pi_0(Y_0)|^2 + \int_0^T \ell (Y_t,V_t,U_t\,;X_t) \ud t \Big) \nonumber\\
&\text{Subj. }\  -\ud Y_t(x) = \big(\clA Y_t(x)+h^\top(x)(U_t + V_t(x))\big) \ud t - V_t^\top(x) \ud Z_t\nonumber \\
&\quad\quad\quad Y_T(x) = f(x) \; \forall x\in \bS \label{eq:opt-cont-duality}
\end{align}
where the set of admissible control ${\cal U} :=
L^2_\clZ([0,T],\Re^m)$ and the cost function 
\begin{equation*}
l(y,v,u\,;x) = \half Q(y\,;x) + \half (u+v(x))^\top R (u+v(x))
\end{equation*}
where $Q(y\,;x)$ is a certain non-negative function. Explicit formulae of
$Q$ for particular examples (finite state-space and It\^o-diffusions) of Markov processes appear in~\citep[Sec. 2]{kim2019bsdeduality}.
It is noted that the constraint is a backward stochastic differential equation (BSDE) with solution $(Y,V):=\{(Y_t,V_t): t\in[0,T]\}\in L^2_{\clZ}([0,T];C(\bS))\times
L^2_{\clZ}([0,T];C^m(\bS))$. The terminal condition $f \in C(\bS)$ is prescribed. 

Consider the following linear structure of the estimator:
$$
S_T = \pi_0(Y_0) - \int_0^T U_t^\top \ud Z_t
$$
where $U\in\clU$ is an admissible control and $Y_0$ is obtained by~\eqref{eq:opt-cont-duality}. The precise duality relationship is as follows:

\medskip

\begin{proposition}[Prop. 1 in~\citep{kim2019bsdeduality}]\label{prop:minimum_variance}
	Consider the observation model~\eqref{eq:obs_model}, together with the dual optimal control
	problem~\eqref{eq:opt-cont-duality}. 
	Then for any choice of admissible control $U\in {\cal U}$: 
	\begin{equation*}
	{\sf J}(U) = \half \E (|S_T - f(X_T)|^2)
	\end{equation*}	
\end{proposition}

\medskip

The significance of the duality relationship is as follows: The problem of obtaining the minimum variance estimate $S_T$
of $f(X_T)$ (minimizer of the right-hand side of the equality) is 
converted into the problem of finding the optimal control $U$
(minimizer of the left-hand side of the identity).  Additional details
including the use of the dual optimal control problem to derive the
nonlinear filter can be found in~\citep{kim2019bsdeduality}.

\section{Proofs of propositions}

\subsection{Proof of Proposition~\ref{thm:duality}}\label{apdx:proof-thm1}

By linearity, $\clL(U;c) = \clL(U;0)+c\ones$ for $U\in\clU$ and $c\in
\Re$.  Therefore, for $\bpi_0\in \clY^\dagger$,
$$
\langle \bpi_0, \clL(U;c) \rangle_\clY =  \langle \bpi_0,\clL(U;0) \rangle_\clY + c\,\bpi_0(\ones) 
$$
Thus, the main calculation is to transform $\langle \bpi_0,\clL(U;0)
\rangle_\clY$. For this purpose, consider~\eqref{eq:NL-ctrl} with
$c=0$ 
and express $\langle \bpi_0,\clL(U;0) \rangle_\clY = \bpi_0(Y_0)$.
Using the It\^{o}-Wentzell formula for measures~\citep[Theorem 1.1]{krylov2011ito},
\begin{align*}
\ud \big(\bpi_t(Y_t)\big) &= \big(\bpi_t(\clA Y_t) \ud t + \bpi_t(h^\top Y_t)\ud Z_t\big) + \bpi_t(h^\top V_t)\ud t \\
&\quad+ \big(\bpi_t(-\clA Y_t - h^\top U_t - h^\top V_t) \ud t + \bpi_t(V_t)\ud Z_t\big)\\
&= -U_t^\top \bpi_t(h)\ud t + \bpi_t(h^\top Y_t + V_t^\top) \ud Z_t
\end{align*}
Integrating both sides,
$$
\bpi_T(Y_T) - \bpi_0(Y_0) = -\int_0^T U_t^\top\bpi_t(h)\ud t + \int_0^T \bpi_t(h^\top Y_t+V_t^\top) \ud Z_t
$$
Under the probability measure $\tilde{\bP}$, $Z$ is a Wiener process. 
Hence, 
\begin{equation*}\label{eq:pf_thm1_1}
\bpi_0(Y_0) = \hE\Big(\int_0^T U_t^\top \bpi_t(h)\ud t \Big) = \langle \bpi(h),U\rangle_\clU
\end{equation*}
Therefore, 
$$
\langle \bpi_0,\clL(U;c) \rangle_\clY =\langle \bpi(h), U \rangle_\clU + c\,\bpi_0(\ones)
$$

\subsection{Proof of Theorem~\ref{prop:equivalent-definition}}\label{apdx:unobservable-direction} 

\eqref{eq:O1} and~\eqref{eq:O2} are equivalent by the closed range theorem.
The proof of \eqref{eq:O2} $\iff$~\eqref{eq:emu-enu} is presented next.

\newP{Necessity} We first show~\eqref{eq:emu-enu}
$\Rightarrow$~\eqref{eq:O2}. For a given $\bpi_0 \in \Nsp(\clL^\dagger)$, then for any $\mu,\nu\in{\cal P}(\bS)$ such that $\varepsilon\bpi_0 = \mu-\nu$ for some constant $\varepsilon\neq 0$, we have:
$$
\varepsilon\bpi_t(h) = \sigma_t^\mu(h)-\sigma_t^\nu(h)
$$
Since $\bpi_0 \in \Nsp(\clL^\dagger)$ implies $\bpi_t(h) \equiv 0$ $\tilde{\sP}$-a.s., 
\begin{equation}\label{eq:unnormalize-equiv}
\sigma_t^\mu(h) = \sigma_t^\nu(h)\quad t\text{-a.e.}\;\;\tilde{\sP}\text{-a.s.}
\end{equation}
Using the Zakai Eq.~\eqref{eq:Zakai} with $f=\ones$ (the constant function),
\begin{equation}\label{eq:Zakai-for-ones}
\sigma_t^\mu(\ones) = 1 + \int_0^t\sigma_s^\mu(h^\top) \ud Z_s
\end{equation}
Therefore~\eqref{eq:unnormalize-equiv} implies that the normalization
constant $\sigma_t^\mu(\ones) = \sigma_t^\nu(\ones)$ for all
$t\in[0,T]$. Thus, using~\eqref{eq:normalize-filter},
$$
\pi_t^\mu(h) = \pi_t^\nu(h) \quad t\text{-a.e.}\;\;\tilde{\sP}\text{-a.s.}
$$
Finally, this is also $\sP^\mu$-a.s. event since $\sP^\mu \ll \tilde{\sP}$.

\newP{Sufficiency} Assume~\eqref{eq:emu-enu} is not true: There exists $\mu\neq\nu \in \clP(\bS)$ such that $\pi_t^\mu\equiv\pi_t^\nu$. We want to show that $\mu-\nu \in \Nsp(\clL^\dagger)$. Equations~\eqref{eq:Zakai-for-ones} and~\eqref{eq:normalize-filter} are combined into:
$$
\sigma_t^\mu(\ones) = 1+\int_0^t\sigma_s^\mu(\ones)\pi_s^\mu(h) \ud Z_s
$$
This implies $\sigma_t^\mu(\ones)\equiv\sigma_t^\nu(\ones)$, and therefore $\sigma_t^\mu(h)\equiv\sigma_t^\nu(h)$ ($\sP^\mu$-a.s.)  by~\eqref{eq:normalize-filter}. Again it is $\tilde{\sP}$-a.s. event by equivalency. It is obvious that $\mu(\ones)-\nu(\ones)=0$, so $\mu-\nu \in \Nsp(\clL^\dagger)$.

\subsection{Proof of Theorem~\ref{thm:rank}}\label{apdx:proof-rank}

For notational ease, we assume $m=1$.  The objective is to show $\clC
= \Rsp(\clL)$.  The proof below is adapted
from~\cite{peng1994backward}.

The definition of $\Nsp(\clL^\dagger)$ is:
$$
\bpi_0 \in \Nsp(\clL^\dagger) \Leftrightarrow \bpi_0(\ones) = 0 \text{
	and } \bpi_t(h) \equiv 0 \quad \forall\;t\in[0,T]
$$
Since $\Nsp(\clL^\dagger)$ is the annihilator of $\Rsp(\clL)$, we have
$\ones,h \in \Rsp(\clL)$. Consider next the 
Zakai equation~\eqref{eq:Zakai} with the initial condition
$\bpi_0\in\Nsp(\clL^\dagger)$ and $f=h$:
$$
\bpi_t(h) = \bpi_0(h) + \int_0^t \bpi_s(\clA h) \ud s + \int_0^t \bpi_s(h^2) \ud Z_s
$$
Since $t$ is arbitrary, the left-hand side is identically zero for all $t\in[0,T]$ if and only if
$$
\bpi_0(h) = 0,\quad \bpi_t(\clA h) \equiv 0,\quad \bpi_t(h^2)\equiv 0 \quad \forall\;t\in[0,T]
$$
and in particular, this implies $\clA h, h^2\in\Rsp(\clL)$. 

The subspace $\clC$ is obtained by continuing to repeat the steps
ad infinitum: If at the conclusion of the  $k^\text{th}$ step, we find
a function $g\in \clC$ such that $\bpi_t(g)\equiv 0$ for all
$t\in[0,T]$.  Then through the use of the Zakai equation,
$$
\bpi_0(g) = 0,\quad \bpi_t(\clA g) \equiv 0,\quad \bpi_t(hg)\equiv 0
\quad \forall\;t\in[0,T]
$$ 
so $\clA g, hg \in \clC$.  By construction, because $
\bpi_0 \in \Nsp(\clL^\dagger)$, $\clC
= \Rsp(\clL)$.

\subsection{Proof of Proposition~\ref{prop:sufficient}}\label{apdx:injective-h}

\newP{Step 1} We first provide the proof for the case when $m=1$. In
this case, $H$ is a column vector and $H_i$ denotes its
$i^{\text{th}}$ element. 
We claim that:
\begin{equation}\label{eq:subset}
\sp\{\ones, \, H, \, H\cdot H, \, \ldots,\, \underbrace{H\cdot H \cdots H}_{(d-1)\text{ times}}\} = \Re^d
\end{equation}
where (as before) the dot denotes the element-wise product.
Assuming that the claim is true, the result easily follows because the
vectors on left-hand side are contained in $\Rsp(\clL)$
(see~\eqref{eq:obs_gram_nl}). It remains to prove the claim. 
For this purpose, express the left-hand side of~\eqref{eq:subset} as
the column space of the following matrix:
\begin{align*}
\begin{pmatrix}
1 & H_1 & H_1^2 & \cdots & H_1^{d-1}\\
1 & H_2 & H_2^2 &\cdots & H_2^{d-1}\\
\vdots &\vdots  & \vdots &\cdots & \vdots\\
1 & H_d & H_d^2 & \cdots & H_d^{d-1}
\end{pmatrix}
\end{align*}
This matrix is easily seen to be full rank by using the Gaussian elimination:
$$
\begin{pmatrix}
1 & H_1 & H_1^2 & \cdots & H_1^{d-1}\\
0 & H_2-H_1 & H_2^2-H_1^2 &\cdots & H_2^{d-1}-H_1^{d-1}\\
\vdots &\vdots  & \vdots &\cdots & \vdots\\
0 & 0 & 0 & \cdots & \prod_{i=1}^{d-1} (H_d-H_i)\\
\end{pmatrix}
$$
The diagonal elements are non-zero because $H_i\neq H_j$.

\newP{Step 2} In the general case, $H$ is a $d \times m$ matrix and
$H_i$ denotes its $i^\text{th}$ row.  We claim that if $H_i\neq H_j$
for all $i\neq j$ then there exists a vector $\tilde{H}$ in the column
span of $H$ such that $\tilde{H}_i\neq \tilde{H}_j$ for all $i\neq
j$. Assuming that the claim is true, the result follows from the $m=1$
case by considering~\eqref{eq:subset} with $\tilde{H}$.
It remains to prove the claim.  Let $\{e_1,\ldots,e_d\}$ denote the
canonical basis in $\Re^d$.  The assumption means $(e_i-e_j)^\top H$
is a non-zero row-vector in $\Re^m$ for all $i\neq j$. 
Therefore, the null-space of $(e_i-e_j)^\top H$ is a
$(m-1)$-dimensional hyperplane in $\Re^m$. Since there are only ${m
	\choose 2}$ such hyperplanes, there must exist a vector $a\in\Re^m$
such that $(e_i-e_j)^\top Ha \neq 0$ for all $i\neq j$. Pick such an
$a$ and define $\tilde{H} := Ha$.

\newP{Step 3} To show the necessity of the injective property when
$A=0$, assume $H_i=H_j$ for some $i\neq j$. Then the corresponding row
is identical, so it cannot be rank $d$.

\end{document}